\documentstyle[11pt]{article}
\newtheorem{theorem}{{\sc Theorem}}
\newcommand{\bt}{\begin{theorem}}
\newcommand{\et}{\end{theorem}}
\newcommand{\newsection}[1]{\setcounter{equation}{0} \setcounter{theorem}{0}
\section{#1}}

\newcommand{\NI}{\noindent}

\newcommand{\bea}{\begin{eqnarray}}
\newcommand{\eea}{\end{eqnarray}}

\def \b #1 {\bf #1}
\newcommand{\IR}{I\!\!R}
\newcommand{\IE}{I\!\!E}
\newcommand{\IC}{I\!\!C}

\newcommand{\IT}{I\!\!T}
\newcommand{\IN}{I\!\!N}
\newcommand{\IZ}{Z\!\!\!Z}

\newcommand{\cla}{{\cal A}}

\newcommand{\clh}{{\cal H}}
\newcommand{\clp}{{\cal P}}

\newcommand{\clb}{{\cal B}}

\newcommand{\clj}{{\cal J}}
\newcommand{\cln}{{\cal N}}

\newcommand{\clm}{{\cal M}}

\newcommand{\al}{\alpha}

\newcommand{\raro}{\rightarrow}

\newcommand{\vsp}{\vskip 1em}

\newcommand{\ul}{\underline}

\def \qed {\hfill \vrule height6pt width 6pt depth 0pt}
\newcommand{\be}{\begin{equation}}
\newcommand{\ee}{\end{equation}}
\newcommand{\ben}{\begin{eqnarray*}}
\newcommand{\een}{\end{eqnarray*}}
\begin{document}
\thispagestyle {empty}
\sloppy

\centerline{\large \bf Jones index of a quantum dynamical semigroup }

\bigskip
\centerline{\bf Anilesh Mohari }
\smallskip
\centerline{\bf S.N.Bose Center for Basic Sciences, }
\centerline{\bf JD Block, Sector-3, Calcutta-98 }
\centerline{\bf E-mail:anilesh@boson.bose.res.in}
\smallskip
\centerline{\bf Abstract}
\bigskip
In this paper we consider a completely positive map $\tau=(\tau_t,t \ge 0)$ with a faithful normal invariant state $\phi$ 
on a type-$II_1$ factor $\cla_0$ and propose an index theory. We achieve this via a more general Kolmogorov's type of 
construction for stationary Markov processes which naturally associate a nested isomorphic von-Neumann algebras. 
In particular this construction generalizes well known Jones construction associated with a sub-factor of type-II$_1$ 
factor.     

\newpage
\newsection{ Introduction:}

\vsp
Let $\tau=(\tau_t,\;t \ge 0)$ be a semigroup of identity preserving completely
positive normal maps [Da,BR] on a von-Neumann algebra $\cla_0$ acting on a
separable Hilbert space $\clh_0$, where either the parameter $t \in \!R_+$, the
set of positive real numbers or $\!Z^+$, the set of positive integers. In case $t \in \!R_+$, i.e.
continuous, we assume that for each $x \in \cla_0$ the map $t \raro \tau_t(x)$ is continuous in
the weak$^*$ topology. Thus variable $t \in \IT_+$ where $\IT$ is either $\IR$ or $\IN$. We assume further
that $(\tau_t)$ admits a normal invariant state $\phi_0$, i.e. $\phi_0 \tau_t = \phi_0 \forall t \ge 0$.

\vsp
As a first step following well known Kolmogorov's construction of stationary 
Markov processes, we employ GNS method to construct a Hilbert space $\clh$ and an increasing  tower of 
isomorphic von-Neumann type$-II$ factors $\{\cla_{[t}: t \in \;\!R \;\mbox{or}\; \!Z\}$ generated by the weak 
Markov process $(\clh,j_t,F_{t]},t \in \!R \;\mbox{or}\; \!Z, \Omega)$ [BP,AM] where $j_t: \cla_0 \raro \cla_{[t}$ 
is an injective homomorphism from $\cla_0$ into $\cla_{[0}$ so that the projection $F_{t]}=j_t(I)$ is the cyclic space 
of $\Omega$ generated by $\{ j_s(x):-\infty < s \le t,\;x \in \cla_0 \}$. The tower of increasing isomorphic von-Neumann 
algebras $\{\cla_{[t},\; t \in \!R \;\mbox{or}\; \!Z\} $ are indeed a type-II$_{\infty}$ 
factor if and only if $\tau$ is not an endomorphism. In any case the projection $j_0(I)$ is a finite projection in $\cla_{[-t}$ for 
all $t \le 0$. In particular we also find an increasing tower of type-II$_1$ factors $\{\clm_t: t \ge 0 \}$ defined by
$\clm_t= j_0(I)\cla_{[-t}j_0(I)$. Thus Jones in-dices $\{ [\clm_t:\clm_s]:0 \le s \le t \}$ are invariance for the Markov 
semigroup $(\cla_0,\tau_t,t \ge 0, \phi_0)$ and further the map $(t,s) \raro [\clm_t:\clm_s]$ is not continuous if the variable 
$(t,s)$ are continuous i.e. if $\tau=(\tau_t: t \in \IR_+)$. In discrete time dynamical system we find an invariance sequence 
$\{[\clm_{n+1}:\clm_n]: n \ge 0 \}$ canonically associated with the canonical conditional expectation on a sub-factor $\clb_0$ 
of a type-II$_1$ factor $\cla_0$ where $\phi_0$ is the unique normal trace on $\cla_0$. However unlike Jones construction we 
have $[\clm_{n+1}:\clm_n]=d^2$ where $d=[\cla_0:\clb_0]$. This shows that our construction in a sense generalizes two step 
Jones construction in discrete time. A detailed study, needs to be done to explore this new invariance, which seems to be an 
interesting problem!

\vsp
\NI {\bf Acknowledgment: } The author takes the opportunity to acknowledge Prof. Luigi Accardi for an invitation to 
visit Centro Vito Volterra, University of Rome, Tor Vergata during the summer 2005. The author further gratefully 
acknowledge Prof. Roberto Longo and Prof. Francesco Fidaleo for valuable discussion which helped the author to realize that 
the tower of type-II$_1$ sub-factors indeed generalizing well known Jones construction.                 

\newsection{ Stationary Markov Processes and Markov shift: }

\bigskip
A family $(\tau_t,\;t \ge 0)$ of one parameter completely positive maps on a $C^*$ algebra or a von-Neumann 
sub-algebra $\cla_0$ is called a {\it quantum dynamical semigroup} if  
$$\tau_0=I,\;\tau_s \circ \tau_t= \tau_{s+t},\;s,t \ge 0$$ 
Moreover if $\tau_t(I)=I,\;t \ge 0$ it is called a {\it Markov }
semigroup. We say a state $\phi_0$ on $\cla_0$ is {\it invariant } for $(\tau_t)$ 
if $\phi_0(\tau_t(x)) = \phi_0(x)\;\forall t \ge 0$. We fix a Markov semigroup 
$(\cla_0,\tau_t,t \ge 0)$ and also a $(\tau_t)-$invariant state $\phi_0$. 

\vsp
In the following we briefly recall [AM] the basic construction of the minimal forward weak Markov 
processes associated with $(\cla_0,\tau_t,\;t \ge 0,\phi_0)$. The construction goes along the 
line of Kolmogorov's construction of stationary Markov processes or Markov shift with a modification 
[Sa,BP] which takes care of the fact that $\cla_0$ need not be a commutative algebra. Here we review 
the construction given in [AM] in order to fix the notations and important properties.

\vsp
We consider the class $\clm$ of $\cla_0$ valued functions
$\ul{x}: \IT \raro \cla_0$ so that $x_r \neq I$ for finitely
many points and equip with the point-wise multiplication
$(\ul{x}\ul{y})_r=x_ry_r$. We define the map $L: (\clm,\clm)
\raro \IC $ by
\be
L(\ul{x},\ul{y}) =
\phi_0(x_{r_n}^*\tau_{r_{n-1}-r_n}(x_{r_{n-1}}^*(.....x_{r_2}^*
\tau_{r_1-r_2}(x_{r_1}^*y_{r_1})y_{r_2})...y_{r_{n-1}})y_{r_n})
\ee
where $\ul{r}=(r_1,r_2,..r_n)\;r_1 \le r_2 \le .. \le r_n$ is
the collection of points in $\IT$ when either $\ul{x}$ or
$\ul{y}$ are not equal to $I$. That this kernel is well defined
follows from our hypothesis that $\tau_t(I)=I, \; t \ge 0$ and
the invariance of the state $\phi_0$ for $(\tau_t).$ The
complete positiveness of $(\tau_t)$ implies that the map $L$ is a
non-negative definite form on $\clm$. Thus there exists a
Hilbert space $\clh$ and a map $\lambda: \clm \raro \clh$ such
that $$<\lambda(\ul{x}),\lambda(\ul{y}) >= L(\ul{x},\ul{y}).$$
Often we will omit the symbol $\lambda$ to simplify our
notations unless more then one such maps are involved.

\vsp
We use the symbol $\Omega$ for the unique element in $\clh$
associated with $x=(x_r=I,\;r \in \IR )$ and $\phi$ for the 
associated vector state $\phi$ on $B(\clh)$ defined by 
$\phi(X)=<\Omega,X\Omega>$.

\vsp
For each $t \in \IR$ we define shift operator $S_t: \clh \raro
\clh$
by the following prescription:
\be
(S_t\ul{x})_r = x_{r+t}
\ee
It is simple to note that $S = (( S_t ,\;t \in \IR))$ is a unitary
group of operators on $\clh$ with $\Omega$ as an invariant element.

\vsp
For any $t \in \IR$ we set
$$\clm_{t]}= \{\ul{x} \in \clm,\; x_r=I\;\forall r > t \} $$
and $F_{t]}$ for the projection onto $\clh_{t]}$, the closed
linear span of $\{\lambda(\clm_{t]})\}$. For any $x \in \cla_0$
and $t \in \IT$ we also set elements $i_t(x),\in \clm$ defined
by $$i_t(x)_r= \left \{ \begin{array}{ll} x ,&\; \mbox{if}\; r=t
\\ I,&\; \mbox{otherwise}\;  \end{array} \right.$$ 

So the map
$V_+: \clh_0 \raro \clh$ defined by $$V_+x=i_0(x)$$ is an
isometry of the GNS space $\{ x:<x,y>_{\phi_0}=\phi_0(x^*y) \}$
into $\clh$ and a simple computation shows that
$<y,V^*_+S_tV_+x>_{\phi_0}=<y,\tau_t(x)>_{\phi_0}$.  Hence
$$P^0_t=V^*_+S_tV_+,\;t \ge 0$$
where $P^0_tx=\tau_t(x)$ is a contractive semigroup of operators on the GNS
space associated with $\phi_0$.

We also note that $i_t(x) \in \clm_{t]}$ and set $\star$-homomorphisms
$j^0_0: \cla_0 \raro \clb(\clh_{0]})$ defined by
$$
j^0_0(x)\ul{y}= i_0(x)\ul{y}
$$
for all $\ul{y} \in \clm_{0]}.$ That it is well defined follows
from (2.1) once we verify that it preserves the inner product whenever
$x$ is an isometry. For any arbitrary element we extend by linearity.
Now we define $j^f_0: \cla \raro \clb(\clh)$
by
\be
j^f_0(x)=j_0^0(x)F_{0]}.
\ee
Thus $j^f_0(x)$ is a realization of $\cla_0$ at time $t=0$ with
$j^f_0(I)=F_{0]}$. Now we use the shift $(S_t)$ to obtain the
process $j^f=(j^f_t: \cla_0 \raro \clb(\clh),\;t \in \IR )$ and
forward filtration $F=(F_{t]},\;t \in \IR)$ defined by the
following prescription:
\be
j^f_t(x)=S_tj^f_0(x)S^*_t\;\;\;F_{t]}=S_tF_{0]}S^*_t,\;\;t \in \IR.
\ee
So it follows by our construction that $ j^f_{r_1}(y_1)j^f_{r_2}(y_2)... 
j^f_{r_n}(y_n) \Omega = \ul{y}$ where $y_r=y_{r_i},\;$ if $r=r_i$
otherwise $I,\;(r_1 \le r_2 \le .. \le r_n)$.  Thus $\Omega$ is
a cyclic vector for the von-Neumann algebra $\cla$ generated by $\{ j^f_r(x), 
\;r \in \IR, x \in \cla_0 \}$. 

\vsp
From (2.4) we also conclude that $S_tXS^*_t \in \cla$ whenever $X \in \cla$ and
thus we can set a family of automorphism $(\al_t)$ on $\cla$
defined by $$\al_t(X)=S_tXS^*_t$$ Since $\Omega$ is an invariant
element for $(S_t)$, $\phi$ is an invariant state for $(\al_t)$.
Now our aim is to show that the reversible system
$(\cla,\al_t,\phi)$ satisfies (1.1) with $j_0$ as defined in
(2.4), for a suitable choice of $\IE_{0]}$. To that end, for any
element $\ul{x} \in \clm$, we verify by the relation
$<\ul{y},F_{t]}\ul{x}=<\ul{y},\ul{x}>$ for all $\ul{y} \in {\cal M}_{t]}$
that $$(F_{t]}\ul{x})_r=
\left \{ \begin{array}{lll} x_r,&\;\mbox{if}\;r < t;\\
\tau_{r_k-t}(...\tau_{r_{n-1}-r_{n-2}}(\tau_{r_n-r_{n-1}}(
x_{r_n})x_{r_{n-1}})...x_t),&\;\mbox{if}\;r=t\\
I,&\;\mbox{if}\;r > t \end{array} \right. $$
where $r_1 \le ..\le r_k \le t \le .. \le r_n$ is the support of
$\ul{x}$. We also claim that
\be
F_{s]}j^f_t(x)F_{s]}=j^f_s(\tau_{t-s}(x))\;\; \forall s \le t.
\ee
For that purpose we choose any two elements $\ul{y},\ul{y'} \in
\lambda({\cal M}_{s]})$ and check the following steps with the aid of
(2.2): $$ <\ul{y},F_{s]}j^f_t(x)F_{s]}\ul{y'}>=
<\ul{y},i_t(x)\ul{y'}>$$
$$=<\ul{y},i_s(\tau_{t-s}(x))\ul{y'})>.$$ Since $\lambda
(M_{s]})$ spans $\clh_{s]}$ it complete the proof of our claim.

\vsp
We also verify that $<z,V^*_+j^f_t(x)V_+y>_{\phi_0}=
\phi_0(z^*\tau_t(x)y)$, hence
\be
V^*_+j^f_t(x)V_+=\tau_t(x),\;\forall t \ge 0.
\ee

\vsp
For any fix $t \in \IT$ let $\cla_{[t}$ be the von-Neumann algebra generated by the family of operators 
$\{j_s(x): t \le s < \infty,\; x \in \cla_0 \}$. We recall that $j_{s+t}(x)=S^*_tj_s(x)S_t,\; t,s \in \!R$ and 
thus $\alpha_t(\cla_{[0}) \subseteq \cla_{[0}$ whenever $t \ge 0$. Hence $(\alpha_t,\; t \ge 0)$ 
is a E$_0$-semigroup on $\cla_{[0}$ with a invariant normal state $\Omega$ and 
\be
j_s(\tau_{t-s}(x))=F_{s]}\alpha_t(j_{t-s}(x))F_{s]} 
\ee 
for all $x \in \cla_0$. We consider the GNS Hilbert space $(\clh_{ \pi_{\phi_0} }, \pi_{\phi_0}(\cla_0),\omega_0)$ 
associated with $(\cla_0,\phi_0)$ and define a Markov semigroup $(\tau_t^{\pi})$ on 
$\pi(\cla_0)$ by $\tau^{\pi}_t(\pi(x))= \pi(\tau_t(x)$. Furthermore we now identify $\clh_{\phi_0}$ as the 
subspace of $\clh$ by the prescription $\pi_{\phi_0}(x)\omega_0 \raro j_0(x)\Omega$. In such a case $\pi(x)$ 
is identified as $j_0(x)$ and aim to verify for any $ t \ge 0$ that
\be
\tau^{\pi}_t(PXP)=P\alpha_t(X)P 
\ee
for all $X \in \cla_{[0}$ where $P$ is the projection from $[\cla_{[0}\Omega]$ onto 
the GNS space $[j^f_0(\cla_0)\Omega]$ which identified with the GNS space associated with $(\cla_0,\phi_0)$. 
It is enough if we verify for typical elements $X=j_{s_1}(x_1)...j_{s_n}(x_n)$ for any $s_1,s_2,...,s_n \ge 0$ 
and $x_i \in \cla_0$ for $1 \le i \le n$ and $n \ge 1$. We use induction on $n \ge 1$. If $X=j_s(x)$ for some 
$s \ge 0$, (2.8) follows from (2.5). Now we assume that (2.8) is true for any element of the form $j_{s_1}(x_1)...j_{s_n}
(x_n)$ for any $s_1,s_2,...,s_n \ge 0$ and $x_i \in \cla_0$ for $1 \le i \le n$. Fix any $s_1,s_2,,s_n,s_{n+1} 
\ge 0$ and consider $X=j_{s_1}(x_1)...j_{s_{n+1}}(x_{n+1})$. Thus 
$P\alpha_t(X)P=j_0(1)j_{s_1+t}(x_1)...j_{s_{n+t}}(x_{n+1})j_0(1)$. If $s_{n+1} \ge s_n$,
we use (2.5) to conclude (2.8) by our induction hypothesis. Now suppose $ s_{n+1} \le s_n$. In such a case if $s_{n-1} 
\le s_{n}$ we appeal once more to (2.5) and induction hypothesis to verify (2.8) for $X$. Thus we are left to consider 
the case where $s_{n+1} \le s_n \le s_{n-1}$ and by repeating this argument we are left to check only the case
where $s_{n+1} \le s_n \le s_{n-1} \le .. \le s_1 $. But $s_1 \ge 0=s_0$ thus we can appeal to (2.5) at the end of the
string and conclude that our claim is true for any typical element $X$ and hence true for all elements in the $*-$ algebra 
generated by these elements of all order. Thus the result follows by von-Neumann density theorem. We also 
note that $P$ is a sub-harmonic projection [Mo1] for $(\alpha_t:t \ge 0)$ i.e. $\alpha_t(P) \ge P$ for all $t \ge 0$ 
and $\alpha_t(P) \uparrow [\cla_{[0}\Omega]$ as $t \uparrow \infty$.   

\vsp
\NI {\bf THEOREM 2.1: } Let $(\cla_0,\tau_t,\phi_0)$ be a Markov semigroup and $\phi_0$ be $(\tau_t)$-invariant state 
on a $C^*$ algebra $\cla_0$. Then the GNS space $[\pi(\cla_0)\Omega]$ associated with $\phi_0$ can be realized as a closed 
subspace of a unique Hilbert space $\clh_{[0}$ up to isomorphism so that the following hold:

\NI (a) There exists a von-Neumann algebra $\cla_{[0}$ acting on $\clh_{[0}$ and a unital $*$-endomorphism 
$(\alpha_t,\; t \ge 0 )$ on $\cla_{[0}$ with a vector state $\phi(X)=<\Omega,X \Omega>$, $\Omega \in \clh_{[0}$
invariant for $(\alpha_t:t \ge 0)$. 

\NI (b) $P \cla_{[0} P$ is isomorphic with $\pi(\cla_0)''$ where $P$ is the projection from $[\cla_{[0}\Omega]$ onto 
$[j^f(\cla_0)\Omega]$; 

\NI (c) $P\alpha_t(X)P=\tau^{\pi}_t(PXP)$ for all $t \ge 0$ and $X \in \cla_{[0}$; 

\NI (d) The closed span generated by the vectors $\{ \alpha_{t_n}(PX_nP)....\alpha_{t_1}(PX_1P)\Omega: 0 \le t_1 \le t_2 \le ..
\le t_k \le ....t_n, X_1,..,X_n \in \cla_{[0}, n \ge 1 \}$ is $\clh_{[0}$. 

\vsp
\NI {\bf PROOF:} The uniqueness up to isomorphism follows from the minimality property (d). \qed

\vsp
Following the literature [Vi,Sa,BhP,Bh] on dilation we say $(\cla_{[0},\alpha_t,\phi)$ is the minimal E$_0$semigroup associated with 
$(\cla_0,\tau_t,\phi_0)$. We have studied extensively asymptotic behavior of the dynamics $(\cla_0,\tau_t,\phi_0)$ in [AM] and 
Kolmogorov's property of the Markov semigroup introduced in [Mo1] was explored to asymptotic behavior of the dynamics 
$(\cla_{[0},\alpha_t,\phi)$. In particular this yields a criteria for the inductive limit state canonically associated with 
$(\cla_{[0},\alpha_t,\phi)$ to be pure. The notion is intimately connected with the notion of a pure $E_0$-semigroup introduced in
[Po,Ar]. For more details we refer to [Mo2].

\newsection{ Dual Markov semigroup and Time Reverse Markov processes: }

\vsp
Now we are more specific and assume that $\cla_0$ is a von-Neumann algebra and each Markov map $(\tau_t)$ is normal and 
for each $x \in \cla_0$ the map $t \raro \tau_t(x)$ is continuous in the weak$^*$ topology. We assume further that 
$\phi_0$ is also faithful. Following [AM2], in the following we briefly recall the time reverse process associated
with the KMS-adjoint ( or Petz-adjoint ) quantum dynamical semigroup
$(\cla,\tilde{\tau}_t,\phi_0)$. 

\vsp
Let $\phi_0$ be a faithful state and
without loss
of generality let also $(\cla_0,\phi_0)$ be in the standard form
$(\cla_0,J,{\cal P},\omega_0)$
[BR] where $\omega_0
\in \clh_0$, a cyclic and separating vector for $\cla_0$, so
that $\phi_0(x)= <\omega_0,x\omega_0>$ and the closer of the
close-able operator $S_0:x\omega_0
\raro x^*\omega_0, S$ possesses a polar decomposition
$S=J\Delta^{1/2}$ with the self-dual positive
cone $\clp$ as the closure of $\{ JxJx\omega_0:x \in \cla_0 \}$ in
$\clh_0$. Tomita's [BR] theorem says that
$\Delta^{it}\cla_0\Delta^{-it}=\cla_0,\;t
\in \IR$ and $J\cla_0J=\cla'_0$, where $\cla'_0$ is the
commutant of $\cla_0$. We define the modular automorphism group
$\sigma=(\sigma_t,\;t \in \IR )$ on $\cla_0$
by
$$\sigma_t(x)=\Delta^{it}x\Delta^{-it}.$$
Furthermore for any normal state $\psi$ on $\cla_0$
there exists a unique vector $\zeta \in {\cal P}$ so that $\psi(x) = <\zeta,x\zeta>$. 
Note that $\clj \pi(x) \clj \pi(y)\Omega=\clj \pi(x) \Delta^{1 \over 2}\pi(y^*)\Omega$
$=\clj\Delta^{1 \over 2}\Delta^{- {1 \over 2}}\pi(x)\Delta^{1 \over 2}\pi(y^*)\Omega$ 
$=\pi(y)\Delta^{1 \over 2}\pi(x^*)\Delta^{-{1 \over 2}}\Omega$. Thus the Tomita's 
map $x \raro \clj \pi(x)\clj$ is an anti-linear $*$-homomorphism representation of $\cla_0$. This
observation leads to a notion called backward weak Markov processes [AM].  

\vsp
To that end we consider the unique Markov semigroup $(\tau'_t)$ on the commutant
$\cla'_0$ of $\cla_0$ so that $\phi(\tau_t(x)y)=\phi(x\tau'_t(y))$ for all
$x \in \cla_0$ and $y \in \cla'_0$.
We define
weak$^*$ continuous
Markov semigroup $(\tilde{\tau}_t)$ on $\cla_0$ by $\tilde{\tau}_t(x)=J\tau'_t(JxJ)J.$
Thus we have the following adjoint relation
\be
\phi_0(\sigma_{1/2}(x)\tau_t(y))=\phi_0(\tilde{\tau}_t(x)\sigma_{-1/2}(y))
\ee
for all $x,y \in \cla_0$, analytic elements for $(\sigma_t)$. One
can as well describe the adjoint semigroup as Hilbert space adjoint of
a one parameter contractive semigroup $(P_t)$ on a
Hilbert space defined by $P_t:\Delta^{1/4}x\omega_0=\Delta^{1/4}\tau_t(x)\omega_0.$ For more details we refer to [Ci].

\vsp
We also note that $i_t(x) \in \clm_{[t}$ and set $\star$ anti-homomorphisms
$j^b_0: \cla_0 \raro \clb(\clh_{[0})$ defined by
$$
j^b_0(x)\ul{y}= \ul{y}i_0(\sigma_{-{i \over 2}}(x^*))
$$
for all $\ul{y} \in \clm_{[0}.$ That it is well defined follows
from (2.1) once we verify by KMS relation that it preserves the inner product whenever
$x$ is an isometry. For any arbitrary element we extend by linearity.
Now we define $j^b_0: \cla \raro \clb(\clh)$
by
\be
j^b_0(x)=j^b_0(x)F_{[0}.
\ee
Thus $j^b_0(x)$ is a realization of $\cla_0$ at time $t=0$ with
$j^b_0(I)=F_{[0}$. Now we use the shift $(S_t)$ to obtain the
process $j^b=(j^b_t: \cla_0 \raro \clb(\clh),\;t \in \IR )$ and
forward filtration $F=(F_{[t},\;t \in \IR)$ defined by the
following prescription:
\be
j^b_t(x)=S_tj^b_0(x)S^*_t\;\;\;F_{[t}=S_tF_{[0}S^*_t,\;\;t \in \IR.
\ee
A simple computation shows for $-\infty < s \le t < \infty$ that  
\be
F_{[t}j^b_s(x)F_{[t}=j^b_t(\tilde{\tau}_{t-s}(x))
\ee
for all $x \in \cla_0$. It also follows by our construction that 
$j^b_{r_1}(y_1)j^b_{r_2}(y_2)...j^b_{r_n}(y_n)\Omega=\sigma_{-{i \over 2}}(\ul{y})$ 
where $y_r=y_{r_i},\;$ if $r=r_i$ otherwise $I,\;(r_1 \ge r_2 \ge .. \ge r_n)$. 
Thus $\Omega$ is a cyclic vector for the von-Neumann algebra $\cla^b$ generated 
by $\{ j^b_r(x),\;r \in \IR, x \in \cla_0 \}''$. We also von-Neumann algebra 
$\cla^b_{t]}$ generated by $\{ j^b_r(x),\;r \le t, x \in \cla_0 \}''$. The following 
theorems say that there is a duality between the forward and backward weak Markov processes. 

\vsp
\NI {\bf THEOREM 3.1: }[AM] We consider the weak Markov processes $(\cla,
\clh,F_{t]},F_{[t},S_t,j^f_t,\;j^b_t\;\;t \in \IR,\; \Omega)$ associated with $(\cla_0,\tau_t,
\;t \ge 0,\;\phi_0)$
and the weak Markov processes $(\tilde{\cla},\tilde{\clh},\tilde{F}_{t]}, \tilde{F}_{[t},\tilde{S}_t,
\;\tilde{j}^f_t,\;\tilde{j}^b_t,\; t \in \IR,\; \tilde{\Omega})$ associated with $(\cla_0,\tilde{\tau}_t,
\; t \ge 0,\;\phi_0)$. There exists an unique anti-unitary operator $U_0:\clh \raro \tilde{\clh}$ so that

\NI (a) $U_0 \Omega = \tilde{\Omega}$;

\NI (b) $U_0 S_t U^*_0 = \tilde{S}_{-t}$ for all $t \in \IR$;

\NI (c) $U_0 j^f_t(x) U^*_0 = \tilde{j}^b_{-t}(x),\;U_0J^b_t(x)U_0=\tilde{j}^f_{-t}(x)$ for all $t \in \IR$;

\NI (d) $U_0F_{t]}U^*_0=\tilde{F}_{[-t},\;\;U_0F_{[t}U^*_0=\tilde{F}_{-t]}$
for all $t \in \IR$;

\vsp
\NI {\bf THEOREM 3.2: } Let $(\cla_0,\tau_t,\phi_0)$ be as in Theorem 3.1 with $\phi_0$ as faithful. Then the commutant of
$\cla_{[t}$ is $\cla^b_{t]}$ for each $t \in \IR$.

\vsp
\NI {\bf PROOF: } It is obvious that $\cla_{[0}$ is a subset of the commutant of $\cla^b_{0]}$. Note also that $F_{[0}$
is an element in $\cla^b_{0]}$ which commutes with all the elements in $\cla_{[0}$. As a first step note that it is good
enough if we show that $F_{[0}(\cla^b_{0]})'F_{[0} = F_{[0}\cla_{[0}F_{[0}$. As for some $X \in (\cla^b_{0]})'$ and
$Y \in \cla_{[0}$ if we have $XF_{[0} = F_{[0}XF_{[0}=F_{[0}YF_{[0}=YF_{[0}$ then we verify that $XZf=YZf$ where $f$
is any vector so that $F_{[0}f=f$ and $Z \in \cla^b_{0]}$ and thus as such vectors are total in $\clh$ we get $X=Y$ ). Thus
all that we need to show that $F_{[0}(\cla^b_{0]})'F_{[0} \subseteq F_{[0}\cla_{[0}F_{[0}$ as inclusion in other direction
is obvious. We will explore in following the relation that
$F_{0]}F_{[0}=F_{[0}F_{0]}=F_{\{ 0 \}}$ i.e. the projection on the fiber at $0$ repeatedly.
A simple proof follows once we use explicit formulas for $F_{0]}$ and $F_{[0}$ given in [Mo1].

\vsp
Now we aim to prove that $F_{[0}\cla_{[0}'F_{[0} \subseteq F_{[0}\cla^b_{0]}F_{[0}$. Let
$X \in F_{[0}\cla_{[0}'F_{[0}$ and verify that $X\Omega=XF_{0]}\Omega = F_{0]}XF_{0]}\Omega =F_{\{0\}}XF_{\{0\}}\Omega \in
[j^b_0(\cla_0)''\Omega]$. On the other-hand we note by Markov property of the backward process
$(j^b_t)$ that $F_{[0}\cla^b_{0]}F_{[0}=j^b(\cla_0)''$. Thus there exists an element $Y \in \cla^b_{0]}$ so that $X\Omega=Y\Omega$.
Hence $XZ\Omega=YZ\Omega$ for all $Z \in \cla_{[0}$ as $Z$ commutes with both $X$ and $Y$. Since $\{ Z\Omega: Z \in \cla_{[0} \}$
spans $F_{[0}$, we get the required inclusion. Since inclusion in the other direction is trivial as
$F_{[0} \in \cla_{[0}'$ we conclude that $F_{[0}\cla_{[0}'F_{[0} = F_{[0}\cla^b_{0]}F_{[0}.$

\vsp
$F_{[0}$ being a projection in $\cla^b_{0]}$ we verify that $F_{[0}(\cla^b_{0]})'F_{[0} \subseteq (F_{[0}\cla^b_{0]}F_{[0})'$
and so we also have $F_{[0}(\cla^b_{0]})'F_{[0} \subseteq ( F_{[0} \cla_{[0}'F_{[0} )' $ as $\cla^b_{0]} \subseteq \cla_{[0}'$.
Thus it is enough if we prove that
$$F_{[0} \cla_{[0}'F_{[0} = (F_{[0} \cla_{[0} F_{[0})'$$
We will verify the non-trivial inclusion for the above equality. Let $X \in (F_{[0} \cla_{[0} F_{[0})'$
then $X\Omega=XF_{0]}\Omega=F_{0]}XF_{0]}\Omega = F_{\{0\}}XF_{\{0\}}\Omega  \in [j^b_0(\cla_0)\Omega]$.
Hence there exists an element $Y \in F_{[0} \cla_{[0}'F_{[0}$ so that $X\Omega=Y\Omega$. Thus for any $Z \in \cla_{[0}$
we have $XZ\Omega=YZ\Omega$ and thus $XF_{[0}=YF_{[0}$. Hence $X=Y \in F_{[0} \cla_{[0}'F_{[0}$. Thus we get the required inclusion.

\vsp
Now for any value of $t \in \IR$ we recall that $\alpha_t(\cla_{[0})=\cla_{[t}$ and $\alpha_t(\cla_{[0})'=
\alpha_t(\cla_{[0}')$, $\alpha_t$ being an automorphism. This completes the proof as $\alpha_t(\cla^b_{0]})=\cla^b_{t]}$
by our construction.  \qed

\newsection{ Subfactors: } 

\vsp
In this section we will investigate further the sequence of von-Neumann algebra $\{ \cla_{[t}:\;\; t \in \IR \}$
defined in the last section with an additional assumption that $\phi_0$ is also faithful and thus
we also have in our hand backward von-Neumann algebras $\{ \cla^b_{t]}: \;\; t \in \IR \}$.

\vsp
\NI {\bf PROPOSITION 4.1: } Let $(\cla_0,\tau_t,\phi_0)$ be a Markov semigroup with a faithful normal invariant state
$\phi_0$. If $\cla_0$ is a factor then $\cla_{[0}$ is a factor. In such a case the following also hold:

\NI (a) $\cla_0$ is type-I (type-II, type-III) if and only if $\cla_{[0}$ is type-I
(type -II , type-III) respectively;

\NI (b) $\clh$ is separable if and only if $\clh_0$ is separable;

\NI (c) If $\clh_0$ is separable then $\cla_0$ is hyper-finite if and only if $\cla_{[0}$ is hyper-finite.

\vsp
\NI {PROOF:} We first show factor property of $\cla_{[0}$. Note that the von-Neumann algebra $\cla^b_{0]}$ generated by the
backward process $\{ j^b_s(x):s \le 0, x \in \cla_0 \}$ is a sub-algebra of $\cla'_{[0}$, the commutant of $\cla_{[0}$. We fix
any $X \in \cla_{[0} \bigcap \cla'_{[0}$ in the center. Then for any $y \in \cla_0$ we verify that
$X j_0(y) \Omega= XF_{0]}j_0(y)\Omega=F_{0]}XF_{0]}j_0(y)\Omega = j_0(xy)\Omega$ for some $x \in \cla_0$.
Since $Xj_0(y)=j_0(y)X$ we also have $j_0(xy)\Omega=j_0(yx)\Omega$. By faithfulness of the state $\phi_0$ we
conclude $xy=yx$ thus $x$ must be a scaler. Thus we have $Xj_0(y)\Omega=c j_0(y)\Omega$ for some scaler $c \in \IC$.
Now we use the property that $X$ commutes with forward process $j_t(x):\;x \in \cla_0, t \ge 0$ and as well as the
backward processes $\{ j^b_t(x),\; t \le 0 \}$ to conclude that $X \lambda(t,x)= c \lambda(t,x)$. Hence $X=c$. Thus
$\cla_{[0}$ is a factor.

\vsp
Now if $\cla_0$ is a type-I factor, then there exists a non-zero minimal projection $p \in \cla_0$. In such a case we claim
that $j_0(p)$ is also a minimal projection in $\cla_{[0}$. To that end let $X$ be any projection in $\cla_{[0}$ so
that $X \le j_0(p)$. Since $F_{0]}\cla_{[0}F_{0]}=j_0(\cla_0)$ we conclude that $F_{0]}XF_{0]}=j_0(x)$ for some $x \in \cla_0$. Hence
$X =j_0(p)Xj_0(p)=F_{0]}Xj_0(p)=j_0(xp)=j_0(px)$
Thus by faithfulness of the state $\phi_0$ we conclude that $px=xp$. Hence $X=j_0(q)$ where $q$ is a projection smaller then equal to $p$.
Since $p$ is a minimal projection in $\cla_0$, $q=p$ or $q=0$ i.e. $X=j_0(p)$ or $0$. So $j_0(p)$ is also a minimal projection.
Hence $\cla_{[0}$ is a type-I factor. For the converse statement we trace the argument in the reverse direction. Let $p$ be a
non-zero projection in $\cla_0$ and claim that there exists a minimal projection $q \in \cla_0$ so that $0 < q \le p$. Now
since $j_0(p)$ is a non-zero projection in a type-I factor $\cla_{[0}$ there exists a non-zero projection $X$ which is minimal in $\cla_{[0}$
so that $0 < X \le j_0(p)$. Now we repeat the argument to conclude that $X=j_0(q)$ for some projection $q$. Since $X \neq 0$ and minimal,
$q \neq 0$ and minimal in $\cla_0$. This completes the proof for type-I case. We will prove now the case for Type-II.

\vsp
Let $\cla_{[0}$ be type-II then there exists a finite projection $X \le F_{0]}.$ Once more $X=F_{0]}XF_{0]}=j_0(x)$ for some projection
$x \in \cla_0$. We claim that $x$ is finite. To that end let $q$ be another projection so that $q \le x$ and $q=uu^*$ and $u^*u=x$. Then
$j_0(q) \le j_0(x)=X$ and $j_0(q)=j_0(u)j_0(u)^*$ and $j_0(x)=j_0(u)^*j_0(u)$. Since $X$ is finite in $\cla_{[0}$ we conclude that
$j_0(q)=j_0(x)$. By faithfulness of $\phi_0$ we conclude that $q=x$, hence $x$ is a finite projection. Since $\cla_0$ is not type-I, it
is type-II. For the converse let $\cla_0$ be type-II. So $\cla_{[0}$ is either type-II or type-III. We will rule out that the possibility
for type-III. Suppose not, i.e. if $\cla_{[0}$ is type-III, for every projection $p \ne 0$, there exists $u \in \cla_{[0}$ so
that $j_0(p)=uu^*$ and $F_{0]}=u^*u$. In such a case $j_0(p)u=uF_{0]}$. Set $j_0(v)=F_{0]}uF_{0]}$ for some $v \in \cla_0$.
Thus $j_0(pv)=j_0(v)$.  Once more by faithfulness of the normal state $\phi_0$, we conclude $pv=v$. So $j_0(v)=uF_{0]}$.
Hence $j_0(v^*v)=F_{0]}$. Hence $v^*v=1$ by faithfulness of $\phi_0$. Since this is true for any non-zero projection $p$ in
$\cla_0$, $\cla_0$ is type-III, which is a contradiction. Now we are left to show the statement for type-III, which is true
since any factor needs to be either of these three types. This completes the proof for (a).

\vsp
(b) is obvious if $\IT$ is $\IZ$. In case $\IT=\IR$, we use our hypothesis that the map $(t,x) \raro \tau_t(x)$ is
sequentially jointly continuous with respect to weak$^*$ topology.

\vsp
For (c) we first recall from [Co] that hyper-finiteness property, being equivalent to injective property of von-Neumann algebra,
is stable under commutant and countable intersection operation when they are acting on a separable Hilbert space.
Let $\cla_0$ be hyper-finite and $\clh_0$ be separable. We will first prove $\cla_{[0}$ is hyper-finite when $\IT=\IZ$, i.e.
time variable are integers. In such a case for each $n \ge 0$, $j_n$ being injective, $j_n(\cla_0)''= \{ j_0(x): x \in \cla_0 \}''$
is a hyper-finite von-Neumann algebra. Thus $\cla_{[0}=\{ j_n(\cla_0)'': n \ge 0 \}''$ is also hyper-finite as they are acting on a
separable Hilbert space. In case $\IT=\IR$, for each $n \ge 1$ we set von-Neumann sub-algebras $\cln^n_{[0} \subseteq \cla_{[0}$
generated by the elements $\{ j_t(\cla_0)'': t = {r \over {2^n}}, 0 \le r \le n2^n \}$. Thus each $\cla^n_{[0}$ is hyper-finite.
Since $\cla_{[0}' = \bigcap_{n \ge 0} (\cla^n_{[0})'$ by weak$^*$ continuity of the map $t \raro \tau_t(x)$, we conclude
that $\cla_{[0}$ is also hyper-finite being generated by a countable family of increasing hyper-finite von-Neumann algebras.

\vsp
For the converse we recall for a factor $\clm$ acting on a Hilbert space $\clh$, Tomiyama's property ( i.e. there
exists a norm one projection $E: \clb(\clh) \raro \clm$, see [BR1] page-151 for details ) is equivalent to hyper-finite property.
For a hyper-finite factor $\cla_{[0}$, $j_0(\cla_0)$ is a factor in the GNS space identified with the subspace $F_{0]}$. Let $E$
be the norm one projection from $\clb(\clh_{[0})$ on $\cla_{[0}$ and verify that the completely positive map $E_0: \clb(\clh_0)
\raro \cla_0$ defined by $E_0(X)= F_{0]}E(F_{0]}XF_{0]})F_{0]}$ is a norm one projection from $\clb(F_{0]})$ to $\cla_0$. This
completes the proof for (b). \qed

\vsp
\NI {\bf PROPOSITION 4.2: } Let $(\cla_0,\tau_t,\phi_0)$ be a dynamical system as in Proposition 4.1. If $\cla_{[0}$ is a
type-II$_1$ factor which admits a unique normalize faithful normal tracial state then the following hold:

\NI (a) $F_{t]}=I$ for all $t \in \IR$;

\NI (b) $\tau=(\tau_t)$ is a semigroup of $*-$endomorphisms.

\NI (c) $\cla_{[0}=j_0(\cla_0)$.

\vsp
\NI {\bf PROOF:} Let $tr_0$ be the unique normalize faithful normal trace on $\cla_{[0}$. For any fix $t \ge 0$ we set a normal
state $\phi_t$ on $\cla_{[0}$ by $\phi_t(x)=tr_0(\alpha_t(x))$. It is simple to check that it is also a faithful normal trace.
Since $\alpha_t(I)=I$, by uniqueness $\phi_t=tr_0$. In particular $tr_0(F_{0]})=tr_0(\alpha_t(F_{0]})=tr_0(F_{t]})$, by faithful
property $F_{t]}=F_{0]}$ for all $t \ge 0$. Since $F_{t]} \uparrow 1$ as $t \raro \infty$ we have $F_{0]}=I$.
Hence $F_{t]}=\alpha_t(F_{0]})=I$ for all $t \in \IR$. This proves (a). For (b) and (c) we recall that
$F_{0]}j_t(x)F_{0]}=j_0(\tau_t(x))$ for all $t \ge 0$ and $j_t:\cla_0 \raro \cla_{[t}$ is an injective $*-$ homomorphism.
Since $F_{t]}=F_{0]}=I$ we have $j_t(x)=F_{0]}j_t(x)F_{0]}= j_0(\tau_t(x))$. Hence $\cla_{[0}=j_0(\cla_0)$
and $j_0(\tau_t(x)\tau_t(y))=j_0(\tau_t(xy))$ for all $x,y \in \cla_0$. Now by injective property of $j_0$,
we verify (b). This completes the proof. \qed

\vsp
We fix a type-II$_1$ factor $\cla_0$ which admits a unique normalize faithful normal tracial state. Since $\cla_{[0}$ is a type-II
factor whenever $\cla_0$ is so, we conclude that $\cla_{[0}$ is a type-II$_\infty$ factor whenever $\tau_t$ is not an endomorphism
on a such a type-II$_1$ factor. The following proposition says much more.

\vsp
\NI {\bf PROPOSITION 4.3:}  Let $\cla_0$ be a type-II$_1$ factor with a unique normalize normal trace and
$(\cla_0,\tau_t,\phi_0)$ be a dynamical system as in Proposition 4.1. Then the following hold:

\NI (a) $j_0(I)$ is a finite projection in $\cla_{[-t}$ for all $t \ge 0$.

\NI (b) For each $t \ge 0$ $\clm_t=j_0(I)\cla_{[-t}j_0(I)$ is a type-II$_1$ factor and $\clm_0 \subseteq \clm_s ...\subseteq \clm_t \subseteq ..,
\;t \ge s \ge 0$ are acting on Hilbert space $F_{0]}$ where $\clm_0=j_0(\cla_0)$.

\vsp
\NI {\bf PROOF:} By Proposition 4.1 $\cla_{[0}$ is a type-II factor. Thus $\cla_{[0}$ is either type-II$_1$ or type-II$_\infty$. In case
it is type-II$_1$, Proposition 4.2 says that $\cla_{[-t}$ is $j_0(\cla_0)$, hence the statements (a) and (b) are true with
$\clm_t = j_0(\cla_0)$. Thus it is good enough if we prove (a) and (b) when $\cla_{[0}$ is indeed a type-II$_\infty$ factor.
To that end for any fix $t \ge 0$ we fix a normal faithful trace $tr$ on $\cla_{[-t}$ and consider the normal map
$x \raro j_0(x)$ and thus a normal trace trace on $\cla_0$ defined by $x \raro tr(j_0(x))$ for $x \in \cla_0$.
It is a normal faithful trace on $\cla_0$ and hence it is a scaler multiple of the unique trace on $\cla_0$.
$\cla_0$ being a type-II$_1$ factor, $j_0(I)$ is a finite projection in $\cla_{[-t}$. Now the general theory
on von-Neumann algebra [Sa] guarantees that $\clm_t$ is type-II$_1$ factor and inclusion follows as $\cla_{-s]}
\subseteq \cla_{-t]}$ whenever $t \ge s$. That $j_0(\cla_0)=j_0(I)\cla_{[0}j_0(I)$ follows from Proposition 4.1. \qed

\vsp
We have now one simple but useful result.

\NI {\bf COROLLARY 4.4: } Let $(\cla_0,\tau_t,\phi_0)$ be as in Proposition 4.1. Then one of the following statements
are false:

\NI (a) $\cla=\clb(\clh)$

\NI (b) $\cla_0$ is a type-II$_1$ factor.

\vsp
\NI {\bf PROOF :} Suppose both (a) and (b) are true. Let $\phi_t$ be the unique normalized trace on $\clm_t$. As
they are acting on the same Hilbert space, we note by uniqueness that $\phi_t$ is an extension of $\phi_s$ for
$t \ge s$. Thus there exists a normal extension of $(\phi_t)$ to weak$^*$ completion $\clm$ of $\bigcup_{t \ge 0} \clm_t$
( here we can use Lemma 13 page 131 [Sc] ). However if $\cla=\clb(\clh)$, $\clm$ is equal to 
$\clb(\clh_{0]})$. $\clh_{0]}$ being an infinite dimensional Hilbert space we arrive at a contradiction. \qed

\vsp
In case $\clm$ in Corollary 4.4 is a type-II$_1$ factor, by uniqueness of the tracial state we claim that 
$\lambda(t)\lambda(s)=\lambda(t+s)$ where $\lambda(t)=tr(F_{-t]})$ for all $t \ge 0$. The claim follows as
von-Neumann algebra $\clm$  is isomorphic with $\alpha_{-t}(\clm)$ which is equal to $F_{-t]}\clm F_{-t]}$.   
The map $t \raro \lambda(t)$ being continuous we get $\lambda(t)=exp(\lambda t)$ for some $\lambda \le 0$. If $\lambda=0$ 
then $\lambda(t)=1$ so by faithful property of the trace we get $F_{-t]}=F_{0]}$ for all $t \ge 0$. Hence we conclude that 
$(\tau_t)$ is a family of endomorphisms by Proposition 4.2. Now for $\lambda < 0$ we have $tr(F_{-t]}) \raro 0$ as $t \raro \infty$. 
As $F_{-t]} \ge |\Omega><\Omega|$ for all $t \ge 0$, we draw a contradiction. Thus the weak$^*$ completion of $\bigcup_{t \ge 0 }\clm_t$ 
is a type-II$_1$ factor if and only if $(\tau_t)$ is a family of endomorphism. In otherwords if $(\tau_t)$ is not a family of endomorphism 
then the weak$^*$ completion of $\bigcup_{t \ge 0 } \clm_t$ is not a type-II$_1$ factor and the tracial state though exists 
on $\clm$ is not unique.

\newsection{ Jones index of a quantum dynamical semigroup on II$_1$ factor: }

\vsp
We first recall Jones's index of a sub-factor originated to understand the structure 
of inclusions of von Neumann factors of type ${\rm II}_1$. Let $N$ be a sub-factor of a finite factor $M$. $M$ acts naturally as 
left multiplication on $L^2(M,tr)$, where $tr$ be the normalize normal trace. The projection $E_0=[N\omega] \in N'$, where $\omega$ 
is the unit trace vector i.e. $tr(x)=<\omega,x\omega>$ for $x \in M$, determines a conditional expectation $E(x)=E_0xE_0$ on $N$. 
If the commutant $N'$ is not a finite factor, we define the index $[M : N]$ to be infinite. In case $N'$ is also a finite factor, 
acting on $L^2(M,tr)$, then the index $[M : N]$ of sub-factors is defined as $tr(E_0)^{-1}$, which is the Murray-von Neumann 
coupling constant [MuN] of $N$ in the standard representation $L^2(M,tr)$. Clearly index is an invariance for the sub-factors. 
Jones proved $[M : N] \in \{4\cos^2(\pi/ n): n=3,4, \cdots\}\cup[4,\infty]$ with all values being realized for some inclusion 
$N \subseteq M$.   

\vsp
In this section we continue our investigation in the general framework of section 4 and study the case when
$\cla_0$ is type-II$_1$ which admits a unique normalize faithful normal tracial state and
$(\tau_t)$ is not an endomorphism on such a type-II$_1$ factor. By Proposition 4.3 $\cla_{[0}$ is a type-II$_\infty$ factor
and $(\clm_t :t \ge 0)$ is a family of increasing type-II$_1$ factor where $\clm_t= j_0(I)\cla_{[-t}j_0(I)$ for all
$t \ge 0$. Before we prove to discrete time dynamics we here briefly discuss continuous case.  Thus the map $I: (t,s) 
\raro [\clm_t:\clm_s],\;\;0 \le s \le t $ is an invariance for the Markov semigroup $(\cla_0,\tau_t,\phi_0)$. By our definition 
$I(t,t)=1$ for all $t \ge 0$ and range of values Jones's index also says that the map $(s,t) \raro I(s,t)$ is not continuous at 
$(s,s)$ for all $s \ge 0$. Being a discontinuous map we also claim the map $(s,t) \raro I(s,t)$ is not time homogeneous, i.e. 
$I(s,t) \neq I(0,t-s)$ for some $ 0 \le s \le t$. If not we could have $I(0,s+t)=I(0,s)I(s,s+t)=I(0,s)I(0,t)$, i.e. $I(0,t) = 
exp(\lambda t)$ for some $\lambda$, this leads to a contradiction. The non-homogenity suggest that $I$ is far from being simple.
We devote rest of the section discussing a much simple example in discrete time dynamics. 
   
\vsp
To that end we review now Jones's construction [Jo, OhP]. Let $\cla_0$ be a type-II$_1$ factor and $\phi_0$ be the unique normalize 
normal trace. The algebra $\cla_0$ acts on $L^2(\cla_0,\phi_0)$ by left multiplication $\pi_0(y)x=yx$ for $x \in L^2(\cla_0,\phi_0)$. 
Let $\omega$ be the cyclic and separating trace vector in $L^2(\cla_0,\phi_0)$. The projection $E_0=[\clb_0\omega]$ induces a trace 
preserving conditional expectation $\tau: a \raro E_0aE_0$ of $\cla_0$ onto $\clb_0$. Thus $E_0\pi_0(y)E_0=E_0\pi_0(E(y))E_0$ for 
all $y \in \cla_0$. Let $\cla_1$ be the von-Neumann algebra $\{ \pi_0(\cla_0),E_0 \}''$. $\cla_1$ is also a type-II$_1$ factor and 
$\cla_0 \subseteq \cla_1$, where we have identified $\pi_0(\cla_0)$ with $\cla_0$. Jones proved that $[\cla_1 : \cla_0] = 
[\cla_0 : \clb_0]$. Now by repeating this canonical method we get an increasing tower of type-II$_1$ factors $\cla_1 \subseteq 
\cla_2 ...$ so that $[\cla_{k+1} : \cla_k]=[\cla_0 : \clb_0]$ for all $k \ge 0$. Thus the natural question: Is Jones tower 
$\cla_0 \subseteq \cla_1 \subseteq ... \subseteq \cla_k ...$ related with the tower $\clm_0 \subseteq \clm_1 ... \clm_k 
\subseteq \clm_{k+1}$ defined in Proposition 4.3 associated with the dynamics $(\cla_0,\tau_n,\phi_0)$?     

\vsp
To that end recall the von-Neumann sub-factors $\clm_0 \subseteq \clm_1$ and the induced representation of $\clm_1$ 
on Hilbert subspace $H_{[-1,0]}$ generated by $\{ j_0(x_0)j_{-1}(x_{-1})\Omega: x_0,x_{-1} \in \cla_0 \}$. $\Omega$ 
is the trace vector for $\clm_0$ i.e. $\phi_0(x)=<\Omega,j_0(x)\Omega>$. However the vector state given by $\Omega$ 
is not the trace vector for $\clm_1$ as $\clm_1 \neq \clm_0$ ( If so we check by trace property that 
$\phi_0(\tau(x)y\tau(z)) = \phi(j_0(x)j_{-1}(y)j_0(z)) = \phi_0(\tau(zx)y)$ for any $x,y,z \in \cla_0$ 
and so $\tau(zx)=\tau(z)\tau(x)$ for all $z,x \in \cla_0$. Hence by Proposition 4.2 we have $\clm_1=
\clm_0$).  Nevertheless $\clm_1$ being a type-II$_1$ factor there exists a unique normalize trace on 
$\clm_1$.   

\vsp
\NI {\bf PROPOSITION 5.1: } $\clm_1 \equiv \cla_2$ and $[\clm_1 : \clm_0]=d^2$ where $d=[\cla_0 : \clb_0]$.  

\vsp
\NI {\bf PROOF: } Let $\phi_1$ be the unique normalize normal trace on $\cla_1$ and $\clh_1 = L^2(\cla_1, \phi_1)$. We 
consider the left action $\pi_1(x):y \raro xy$ of $\cla_1$ on $\clh_1$. Thus $\pi_0(\cla_0)$ is also acting on $\clh_1$. 
Since $E_0\pi_0(x)E_0 = E_0\pi_0(\tau(x))E_0=E_0\pi_0(\tau(x))$, for any element $X \in \cla_1$, $E_0X=E_0\pi_0(x)$ for 
some $x \in \cla_0$. Thus $\pi_1(E_0)$ is the projection on the subspace $\{ E_0\pi_0(x): x \in \cla_0 \}$.  

\vsp
For any $y \in \cla_0$ we set

\NI (a) $k_{-1}(y)$ on the subspace $\pi_1(E_0)$ by $k_{-1}(y)E_0\pi_0(x)=E_0\pi_0(yx)$ for $x \in \cla_0$ and extend it to 
$\clh_1$ trivially. That $k_{-1}(y)$ is well defined and an isometry for an isometry $y$ follows from the following 
identities: 
$$\phi_1((E_0\pi_0(yz))^*E_0\pi_0(yx)) =\phi_1(\pi_0(z^*y^*)E_0\pi_0(yx))$$ 
$$=\phi_1(E_0\pi_0(yx)\pi_0(z^*y^*))\;\; \mbox{by trace property}$$
$$ =\phi_1(E_0)\phi_0(\pi_0(yx)\pi_0(z^*y^*))\;\;\mbox{being a trace and } \phi_1(E_0\pi_0(x))=
\phi_1(E_0)\phi_0(\pi_0(x))$$ 
$$=\phi_1(E_0)\phi_0(\pi_0(z^*y^*)\pi_0(yx)) =\phi_1(E_0)\phi_0(\pi_0(z^*x)) $$
$$=\phi_1((E_0\pi_0(z))^*E_0\pi_0(x))$$   

\NI (b) $k_0(y)x = \pi_0(y)x$ for $x \in \cla_1$. Thus $y \raro k_0(y)$ is an injective 
$*$-representation of $\cla_0$ in $L^2(\cla_1,\phi_1)$.    

\vsp
For $y,z \in \cla_0$ we verify that 
$$<E_0\pi_0(y), k_{-1}(1)k_0(x)k_{-1}(1)E_0\pi_0(z)>_1 = <E_0\pi_0(y),E_0\pi_0(x)E_0\pi_0(z)>_1 $$ 
$$= <E_0\pi_0(y),E_0\pi_0(\tau(x))E_0\pi_0(z)>_1$$
$$= <E_0\pi_0(y),E_0\pi_0(\tau(x))\pi_0(z)>_1$$
Thus $k_{-1}(1)k_0(x)k_{-1}(1)=k_{-1}(\tau(x))$ for all $x \in \cla_0$. Note that $k_{-1}(1)=\pi_1(E_0)$ and the 
identity operator in $\clh_1$ is a cyclic vector for the von-Neumann algebra $\{k_0(x),k_{-1}(x),x \in \cla_0 \}''$. 
We have noted before that the vector $\Omega$ need not be the tracial vector for $\clm_1$ and also verify by a direct 
computation that the space $\{k_0(y)k_{-1}(x)1: x,y \in \cla_0 \}$ is equal to $\{yE_0x: y,x \in \cla_0 \}$ which is 
a proper subspace of $L^2(\cla_1,\phi_1)$.

\vsp
Now we claim that the type-II$_1$ factor $\clm_1$ is isomorphic to the von-Neumann algebra 
$\{k_{-1}(x), k_0(x),\;x \in \cla_0) \}''$. To that end we define an unitary operator from 
$L^2(\cla_1,\phi_1)$ to $L^2(\clm_1,tr_1)$ by taking an element $k_{t_1}(x_1)..k_{t_n}(x_n)$ to 
$j_{t_1}(x_1)..j_{t_n}(x_n)$, where $t_k$ are either $0$ or $-1$. That it is an unitary operator 
follows by the tracial property of the respective states and weak Markov property of the homomorphisms. 
We leave the details and without lose of generality we identify these two weak Markov processes. Since 
$k_{-1}(1)=\pi_1(E_0)$, we conclude that $\pi_1(\cla_1) \subseteq \clm_1$. In fact strict inclusion hold 
unless $\clb_0=\cla_0$. 

\vsp
However by our construction $\cla_1=\pi_0(\cla_0) \vee E_0$ is acting on $L^2(\cla_0,\phi_0)$ and $\cla_2 = 
\pi_1(\cla_1) \vee E_1$ is acting on $L^2(\cla_1,\phi_1)$ where $E_1$ is the cyclic subspace of $1$ generated 
by $\pi_1(\pi_0(\cla_0))$ i.e. $E_1=[\pi_1(\pi_0(\cla_0))1]$. From (a) we also have 
\be
k_{-1}(y)\pi_1(E_0)E_1=\pi_1(E_0)k_0(y)E_1
\ee
for all $y \in \cla_0$. 

By Temperley-Lieb relation [Jo] we have $\pi_1(E_0)E_1\pi_1(E_0)={1 \over d}\pi_1(E_0)$ and thus
post multiplying (5.1) by $\pi_1(E_0)$ we have 
\be
{1 \over d } k_{-1}(y)\pi_1(E_0) = \pi_1(E_0)k_0(y)E_1\pi_1(E_0)
\ee 
So it is clear now that $k_{-1}(y) \in \pi_1(\cla_1) \vee E_1$ for all $y \in \cla_0$. Thus $\pi_1(\cla_1) \subseteq 
\clm_1 \subseteq \pi_1(\cla_1) \vee E_1=\cla_2$. 

\vsp 
We claim also that $E_1 \in \clm_1$. We will show that any unitary element $u$ commuting with $\clm_1$ is 
also commuting with $E_1$. By (5.2) we have $\pi_1(E_0)k_0(y)(uE_1u^*-E_1)\pi_1(E_0)=0$ for all $y \in 
\cla_0$. By taking adjoint we have $\pi_1(E_0)(uE_1u^*-E_1)k_0(y)\pi_1(E_0)=0$ for all $y \in \cla_0$.
Since $\cla_1=\pi_0(\cla_0)\vee E_0$ we conclude by cylicity of the trace vector that 
$\pi_1(E_0)(uE_1u^*-E_1)=0$. So we have $E_1\pi_1(E_0)uE_1u^*=E_1\pi_1(E_0)E_1=E_1$ by 
Temperley-Lieb relation. So $u^*E_1u \pi_1(E_0) = u^*E_1\pi_1(E_0)uE_1u^* u = u^*E_1u$
By taking adjoint we get $\pi_1(E_0) u^*E_1u = u^*E_1u$. Since same is true for $u^*$, we conclude
that $u^*E_1u=E_1$ for any unitary $u \in \clm_1'$. Hence $E_1 \in \clm_1$.   

\vsp
Hence $\clm_1=\cla_2$. Since $[\clm_1 : \clm_0]= [\clm_1 : \cla_1] [\cla_1 : \clm_0]$ and 
$[\cla_1 : \cla_0] = [\cla_0 : \clb_0]=d$, we conclude the result. \qed   

\vsp
\NI {\bf THEOREM 5.2: } $\clm_m \equiv \cla_{2m}$ for all $m \ge 1$. 

\vsp
\NI {\bf PROOF: } Proposition 5.1 gives a proof for $m=1$. The proof essentially follows the same steps as in 
Proposition 5.1. We use induction method for $m \ge 1$. Assume it is true for $1,2,..m$. Now consider the Hilbert space 
$L^2(\cla_{2m+1},tr_{2m+1}),\;m \ge 1$ and we set homomorphism $k_{-1},k_0$ from $\cla_{2m}$ into $\clb(L^2(\cla_{2m+1},tr_{2m+1}))$ 
in the following: 

\NI (a) $k_0(x)y=xy$ for all $y \in \cla_{2m+1}$ and $x \in \cla_{2m}$

\NI (b) $\pi_{2m+1}(E_{2m})$ is the projection on the subspace $\{E_{2m}y : y \in \cla_{2m} \}$ 
and $k_{-1}(x)$ defined on the subspace $E_{2m}$ by 
$$k_{-1}(x)E_{2m}y=E_{2m}xy $$ 
for all $x \in \cla_{2m}$ and $y \in \cla_{2m}$. That $k_{-1}$ is an homomorphism follows as in Proposition 5.1. Thus 
an easy adaptation of Proposition 5.1 says that $\clm = \{k_0(x),k_{-1}(x): x \in \cla_{2m} \}''$ is a type-II$_1$ factor and
proof will be complete once we show that it is isomorphic to $\clm_{m+1}$. 

\vsp
To that end we check as in Proposition 5.1 that 
$$k_{-1}(E_{2m})=k_{-1}(I),\;\;k_{-1}(I)k_0(x)k_{-1}(I)=k_{-1}(E_{2m}xE_{2m})$$ 
for all $x \in \cla_{2m}$ and $k_{-1}(x)E_{2m}=E_{2m}k_0(x)E_{2m+1}$ where $E_{2m+1}$ is the cyclic space of the trace vector generated 
by $\cla_{2m}$. Thus following Proposition 5.1 we verify now that type-II$_1$ factor $\clm=\{ k_0(x), k_{-1}(I), E_{2m+1}, x \in \cla_{2m} \}''$ 
is isomorphic to $\clm_{m+1} = \{ J_{-1}(x), J_0(x)\;\;\;x \in \cla_m \}''$, where we used notation $J_{-1}(x)=x$ for all $x \in 
\cla_m,\;\;J_{0}(x)=S J_{-1}(x)S^*$ for all $x \in \cla_{2m}$ where we have identified $\cla_{2m} \equiv \clm_m$ with 
$\{j_k(x),-m-1 \le k \le -1, x \in \cla_0 \}''$ and $S$ is the (right) Markov shift. This completes the proof. \qed

\bigskip
{\centerline {\bf REFERENCES}}

\begin{itemize} 

\bigskip
\item{[AM]} Accardi, L., Mohari, A.: Time reflected Markov processes. Infin. Dimens. Anal. Quantum Probab. Relat. Top., 
vol-2, no-3, 397-425 (1999).

\item {[Ar]} Arveson, W.: Pure $E_0$-semigroups and absorbing states, Comm. Math. Phys. 187 , no.1, 19-43, (1997)

\item {[Bh]} Bhat, B.V.R.: An index theory for quantum dynamical semigroups, Trans. Amer. Maths. Soc. vol-348, no-2 561-583 (1996).   

\item {[BP]} Bhat, B.V.R., Parthasarathy, K.R.: Kolmogorov's existence theorem for Markov processes on $C^*$-algebras, Proc.
Indian Acad. Sci. 104,1994, p-253-262.

\item {[BR]} Bratelli, Ola., Robinson, D.W. : Operator algebras and quantum statistical mechanics, I,II, Springer 1981.

\item{[Da]} Davies, E.B.: Quantum Theory of open systems, Academic press, 1976.

\item{[El]} Elliot. G. A.: On approximately finite dimensional von-Neumann algebras I and II, Math. Scand. 39 (1976), 91-101;
Canad. Math. Bull. 21 (1978), no. 4, 415--418. 

\item{[Jo]} Jones, V. F. R.: Index for subfactors. Invent. Math. 72 (1983), no. 1, 1--25.

\item{[Mo1]} Mohari, A.: Markov shift in non-commutative probability, Jour. Func. Anal. 199 (2003) 189-209.  

\item{[Mo2]} Mohari, A.: Pure inductive limit state and Kolmogorov's property, ....
 
\item{[Mo3]} Mohari, A.: SU(2) symmetry breaking in quantum spin chain, ....

\item{[MuN]} Murray, F. J.; von Neumann, J., On rings of operators. (English)[J] Ann. Math., Princeton, (2)37, 116-229.

\item{[OP]} Ohya, M., Petz, D.: Quantum entropy and its use, Text and monograph in physics, Springer-Verlag 1995. 

\item{[Po]} Powers, Robert T.: An index theory for semigroups of $*$-endomorphisms of
$\clb(\clh)$ and type II$_1$ factors.  Canad. J. Math. 40 (1988), no. 1, 86--114.

\item{[Sak]} Sakai, S.: C$^*$-algebras and W$^*$-algebras, Springer 1971.  

\item{[Sa]} Sauvageot, Jean-Luc: Markov quantum semigroups admit covariant Markov $C^*$-dilations. Comm. Math. Phys. 
106 (1986), no. 1, 91­103.

\item{[Vi]} Vincent-Smith, G. F.: Dilation of a dissipative quantum dynamical system to a quantum Markov process. Proc. 
London Math. Soc. (3) 49 (1984), no. 1, 58­72. 

\end{itemize}

\end{document}